\documentclass[twoside,notitlepage,12pt]{article}

\usepackage{amssymb}
\usepackage[leqno]{amsmath}
\usepackage{amsfonts}
\usepackage{amsopn}
\usepackage{amstext}
\usepackage{amsthm}

\usepackage[colorlinks, backref]{hyperref}

\newcommand{\nice}{\ensuremath{{\mathbb L}_1}}

\providecommand{\cal}{\mathcal}
\renewcommand{\Bbb}{\mathbb}

\newenvironment{pf}{\begin{proof}}{\end{proof}}

\newcommand{\Dee}{{\cal{D}}}
\newcommand{\Eee}{{\cal{E}}}
\newcommand{\Ef}{{\cal{F}}}
\newcommand{\Gee}{{\cal{G}}}

\newcommand{\Pee}{{\cal{P}}}

\newcommand{\Ve}{{\Bbb{V}}}

\newcommand{\Err}{{\Bbb{R}}}

\newcommand{\lam}{{\lambda}}
\newcommand{\al}{\alpha}

\newcommand{\sig}{\sigma}
\newcommand{\eps}{\varepsilon}
\renewcommand{\phi}{\varphi}
\renewcommand{\rho}{\varrho}

\newcommand{\rest}{\restriction}

\newcommand{\ntr}{n\in\omega}

\newcommand{\loe}{\leqslant}
\newcommand{\goe}{\geqslant}

\newcommand{\subs}{\subseteq}

\newcommand{\cl}{\operatorname{cl}}

\newcommand{\id}[1]{{\operatorname{id}_{#1}}} 

\newcommand{\dom}{\operatorname{dom}}

\newcommand{\poset}{{\Bbb{P}}}

\newcommand{\concat}{{}^\smallfrown}

\newcommand{\by}{/}

\newcommand{\borel}[1]{\operatorname{Borel}\left({#1}\right)}

\newtheorem{tw}{Theorem}[section]

\newtheorem{lm}[tw]{Lemma}

\theoremstyle{definition}

\newtheorem{pyt}{Question}

\theoremstyle{remark}
\newtheorem{case}{Case}

\newcommand{\setof}[2]{\{#1\colon #2\}}

\newcommand{\seq}[1]{\langle #1 \rangle}

\newcommand{\seqq}[2]{{\langle #1 \rangle}_{#2}}
\newcommand{\sett}[2]{\{#1\}_{#2}}
\newcommand{\sn}[1]{\{#1\}} 
\newcommand{\dn}[2]{\{#1,#2\}} 
\newcommand{\pair}[2]{\langle #1, #2 \rangle} 
\newcommand{\map}[3]{#1\colon #2 \to #3} 
\newcommand{\img}[2]{#1[#2]} 

\newcommand{\dpower}[2]{[#1]^{#2}}

\newcommand{\fin}[1]{[#1]^{<\omega}}

\providecommand{\nat}{\omega}

\newcommand{\ciag}[1]{{\sett{{#1}_n}{\ntr}}}

\newcommand{\cmp}{\circ}
\newcommand{\Cantor}[1]{2^{#1}}
\newcommand{\cantor}{\Cantor\omega}

\title{Covering an uncountable square by countably many continuous functions}
\author{
{\sc Wies{\l}aw Kubi\'s}\\ \\
{\small Institute of Mathematics}\\
{\small Czech Academy of Sciences, Prague}\\
{\small Czech Republic}\\
{\small\href{mailto:kubis@math.cas.cz}{kubis@math.cas.cz}}
\and
{\sc Benjamin Vejnar}\\ \\
{\small Department of Mathematical Analysis}\\
{\small Charles University, Prague}\\
{\small Czech Republic}
}

\begin{document}

\maketitle

\begin{abstract}
We prove that there exists a countable family of continuous real functions whose graphs together with their inverses cover an uncountable square, i.e. a set of the form $X\times X$, where $X\subs\Err$ is uncountable. This extends Sierpi\'nski's theorem from 1919, saying that $S\times S$ can be covered by countably many graphs of functions and inverses of functions if and only if $|S|\loe\aleph_1$. 
Our result is also motivated by Shelah's study of planar Borel sets without perfect rectangles.

\vspace{3mm}


\end{abstract}

\section{Introduction}

A classical result of Sierpi\'nski from 1919 (see \cite{Sie1919, Sie1924} or \cite[Chapter I]{Sie}) says that, given a set $S$ of cardinality $\aleph_1$, there exists a countable family of functions $\map {f_n}SS$ such that
\begin{equation}
S\times S=\bigcup_{\ntr}(f_n\cup f_n^{-1}),	
\label{eq_one}
\end{equation}
where $f_n^{-1}$ is the inverse of $f_n$, i.e. $f_n^{-1}=\setof{\pair{f_n(x)}{x}}{x\in S}$.
A typical proof proceeds as follows. Assume $S=\omega_1$ and for each positive $\beta<\omega_1$ choose a surjection $\map{g_\beta}{\omega}\beta$. Define $\map{f_n}{\omega_1}{\omega_1}$ by the equation
$f_n(\beta) = g_\beta(n)$.
For every $\pair\al\beta\in S\times S$ with $\al<\beta$ there exists $n$ such that $g_\beta(n)=\al$; thus $\pair \al\beta\in f_n$ and $\pair \beta\al\in f_n^{-1}$. Finally, it suffices to add the identity function to the family $\sett{f_n}{\ntr}$ in order to get (\ref{eq_one}). It is worth noting that the sets $f_n^{-1}(\al)$ form an {\em Ulam matrix} on $\omega_1$. See e.g. \cite[Chapter 10]{Jech} or \cite[Chapter II, \S6]{Kunen} for applications of Ulam matrices.

An easy argument (also noted by Sierpi\'nski) shows that the above statement fails for a set $S$ of cardinality $\aleph_2$. In particular, the continuum hypothesis is equivalent to the statement ``there exists a countable family of functions whose graphs together with their inverses cover the plane".

Let us say that a set $M$ is {\em covered by} a family of functions $\Ef$, if for every $\pair xy\in M$ there is $f\in\Ef$ such that either $y=f(x)$ or $x=f(y)$.

\begin{pyt}\label{pytanko}
Does there exist a sequence $\sett{\map{f_n}\Err\Err}{\ntr}$ of continuous functions that covers an uncountable square?
\end{pyt}


One can hope for a positive answer only when the side $S$ of the square has some smallness properties, besides having cardinality $\aleph_1$.
In fact, by a result of Zakrzewski \cite[Theorem 2.1]{Zak}, if $S\times S$ is covered by countably many functions (and their inverses) whose graphs are Borel sets, then $S$ is {\em universally small}, i.e. $S$ belongs to every Borel $\sig$-ideal $I\subs\Pee(\Err)$ such that $\borel\Err\by I$ satisfies the countable chain condition.

Actually, consistent affirmative answers to Question \ref{pytanko} already exist in the literature. 
Namely, Abraham and Geschke \cite{AG} showed that for every set $X \subs \Err$ of cardinality $\aleph_1$ there
is a ccc forcing notion adding countably continuous functions that cover $X \times X$.
Conseuqently, under Martin's Axiom every $\aleph_1$-square in the plane is covered by a countable family of continuous functions.
The Open Coloring Axiom of Abraham, Rubin and
Shelah \cite{ARS} implies that for every set $X \subs \cantor$ of size $\aleph_1$ there is a countable
family of 1-Lipschitz functions that covers $X \times X$ (see \cite{G} for more details).


Yet another motivation for addressing Question~\ref{pytanko}
comes from the work of Shelah \cite{Sh522}, continued in \cite{KS}, where planar Borel sets without perfect squares were studied. 
It is not hard to prove (see e.g. \cite[Thm. 2.2]{Kubis1}) that a $G_\delta$ subset of the plane containing countable squares of arbitrarily large countable Cantor-Bendixson ranks, contains also a perfect square. On the other hand, using Keisler's absoluteness,
it has been proved in \cite{Sh522} that there exists in ZFC a planar $F_\sig$ set $C$ such that $S\times S\subs C$ for some uncountable set $S$, while $P_0\times P_1\not\subs C$ whenever $P_0$, $P_1$ are perfect sets. A significant part of \cite{KS} is devoted to a ZFC construction of certain $F_\sig$ sets in the plane which do not contain perfect squares, while consistently they contain squares of a prescribed cardinality below $\aleph_{\omega_1}$. 
These sets moreover have certain universality property, among sets of the same type (see \cite{KS} for details). Based on the results of \cite{Sh522} and \cite{KS}, it is natural to ask for the existence of a more special planar $F_\sig$ set which covers an uncountable square: namely, a set consisting of countably many continuous real functions and their inverses. There are natural restrictions here. Namely, such a set cannot contain rectangles of the form $S_0\times S_1$, where 
$|S_i|=\aleph_1$ and $|S_{1-i}|\goe\aleph_2$.
Easy absoluteness arguments show that it cannot contain perfect rectangles, therefore the best property we can expect is covering a square of cardinality $\aleph_1$. 

In the present note we find a family of continuous functions $\Ef=\sett{\map {f_n}{\cantor}{\cantor}}{\ntr}$ such that every maximal square covered by $\Ef$ is uncountable.
The functions $f_n$ are not Lipschitz with respect to any natural metric on $\cantor$, however we describe a natural ccc forcing notion which introduces a countable family of $1$-Lipschitz functions on the Cantor set that covers an uncountable square. Using Keisler's completeness \cite{Kei}, we deduce that such a family exists in ZFC, although we do not know any direct construction.

Finally, we observe that it is impossible to cover the square of any uncountable compact Hausdorff space by countably many continuous functions and their inverses.

It is worth noting that there is {\em no} uncountable $S\subs \Err$ whose square can be covered by countably many non-decreasing functions and their inverses. This is because the graphs of such functions (and of their inverses) are chains with respect to the coordinatewise ordering and consequently the order of $S$ would be a Countryman type (see \cite{Sh50} or \cite[p. 258]{Todorcevic_handbook}), not embeddable into the real line.

\section{Main result}

\begin{tw}\label{wuiehfoa}
There exists a family of continuous functions $\Ef=\sett{\map{f_n}{\cantor}{\cantor}}{\ntr}$ such that every maximal square covered by $\Ef$ is uncountable.
\end{tw}

\begin{pf}
Let $\omega=\bigcup_{\ntr}A_n$, where the sets $A_n$ are infinite and pairwise disjoint. For each $n$ let $\map{\phi_n}\omega{A_n}$ be a bijection. This induces a homeomorphism $\map {h_n}{\Cantor{A_n}}{\cantor}$ given by $h_n(x)(i)=x(\phi_n(i))$, $i<\nat$.
Let $\map {p_n}{\cantor}{\Cantor{A_n}}$ be the canonical projection, i.e. $p_n(x)=x\rest{A_n}$ for $x\in\cantor$. For each $n>0$ define
$$f_n = h_n\cmp p_n.$$
Finally, let $f_0$ be the identity. We claim that the family $\Ef=\sett{f_n}{\ntr}$ has the desired property. 

Let $S$ be a maximal, with respect to inclusion, set whose square is covered by $\Ef$. Suppose $|S|\loe\aleph_0$ and write $S=\sett{x_n}{n>0}$. Consider $y\in\cantor$ satisfying $$y\rest A_n = h_n^{-1}(x_n)\qquad \text{for every }n>0.$$ Notice that there are $2^{\aleph_0}$ possibilities for defining $y\rest A_0$, therefore we may choose $y$ so that it does not belong to $S$. Finally, observe that given $n>0$ we have that $f_n(y) = h_n(y\rest A_n) = x_n$.
This shows that the square of $S\cup\sn y$ is covered by $\Ef$, a contradiction.
\end{pf}

Notice that the Axiom of Choice was heavily used in the above simple proof. The above result says in particular that $\Ef$ covers an uncountable square. One can ask whether the latter statement is a theorem of ZF. The construction of the family $\Ef$ of course does not require the AC. An uncountable set whose square is covered by $\Ef$ can be constructed by transfinite induction and this again does not require the AC. The only problem is that one needs to have a sequence of surjections $\seqq{\map {g_\beta}\omega\beta}{0<\beta<\omega_1}$. The existence of such a sequence clearly implies that $\omega_1$ is a regular cardinal, which is not provable in ZFC.
In fact, the proof of Theorem \ref{wuiehfoa} shows that (using AC) there exists a one-to-one sequence $\seqq{x_\al}{\al<\omega_1}$ in $\cantor$ such that for every $\al<\beta<\omega_1$ there is $n>0$ satisfying $x_\al=f_n(x_\beta)$. 
This property of $S=\sett{x_\al}{\al<\omega_1}$ is formally stronger than saying ``$S^2$ is covered by $\Ef$" and it clearly implies the existence of the above mentioned sequence of surjections $\map{g_\beta}\omega\beta$. Without using the AC, we are only able to conclude the following.

\begin{tw} Assume ZF. There exists a countable family of continuous self-maps of the Cantor set $\cantor$ such that for every countable set $S\subs \cantor$ whose square is covered by $\Ef$ there exists $y\in \cantor\setminus S$ such that the square of $S\cup\sn y$ is also covered by $\Ef$. More precisely, for every $x\in S$ there is $f\in\Ef$ such that $x=f(y)$.
\end{tw}

\begin{pf}
Note that the sentence ``$S$ is countable" means that there exists a bijection $\map x\omega S$.
Repeat the proof of Theorem \ref{wuiehfoa}, defining $y$ specifically, by fixing a bijection $\map \phi\omega {A_0}$ and setting $y(\phi(n))=x_n(\phi(n))+1$. Then $y\ne x_n$ for any $\ntr$.
\end{pf}

\section{Forcing countably many Lipschitz functions}

It is easy to see that the functions constructed in Theorem \ref{wuiehfoa} are not Lipschitz with respect to the natural metric on the Cantor set. We do not know a direct construction of a countable family of $1$-Lipschitz functions covering an uncountable square. This section is devoted to showing that such a family can be introduced by a natural forcing notion. Using Keisler's completeness, we later conclude that this family exists in ZFC, namely:

\begin{tw}\label{wfojwqeoj}
There exists a countable family of $1$-Lipschitz functions on the Cantor set which covers an uncountable square.
\end{tw}

We first show the consistency of the above statement with the axioms of ZFC. The metric on $\cantor$ which we have in mind is given by the formula $d(x,y)=2^{-k}$, where $k$ is the smallest natural number such that $x\rest k\ne y\rest k$.

Given a natural number $n$, we shall denote by $2^n$ the complete binary tree consisting of all zero-one sequences of length $n$. Trees of the form $2^n$ serve as finite approximations of the Cantor set $\cantor$. We consider $2^n$ with the lexicographic ordering and with the metric defined above, like in the case of $\cantor$. Denote by $\nice(n)$ the set of all $1$-Lipschitz functions of the form $\map g{2^n}{2^n}$. 

We are going to define a forcing notion $\poset$ which will introduce a countable family of Lipschitz functions covering an uncountable square. 

A condition $p\in\poset$ is, by definition, of the form $p=\seq{n^p,s^p,v^p,\Ef^p,\gamma^p,\rho^p}$, where 
\begin{enumerate}
	\item[(1)] $n^p\in\nat$, $s^p\in\fin \nat$ and $v^p\in\fin{\omega_1}$;
	\item[(2)] $\Ef^p=\sett{f^p_i}{i\in s^p}\subs{\nice(n^p)}$ and $\map{\rho^p}{\dpower{v^p}2}{s^p}$;
	\item[(2')]$\rho^p(\al,\beta)\ne\rho^p(\al',\beta)$ whenever $\al<\al'<\beta$;
	\item[(3)] $\map{\gamma^p}{v^p}{2^{n^p}}$ is one-to-one;
	\item[(4)] $\gamma^p(\al)=f^p_{\rho^p(\al,\beta)}(\gamma^p(\beta))$ whenever $\al<\beta$ and $\al,\beta\in v^p$.	
\end{enumerate}
Note that condition (2') is actually implied by the conjunction of (3) and (4).
The order of $\poset$ is defined naturally. Namely, $p\loe q$ ($q$ is stronger than $p$) iff
\begin{enumerate}
	\item[(5)] $n^p\loe n^q$, $s^p\subs s^q$, $v^p\subs v^q$;
	\item[(6)] $f^q_i(\eta)\rest n^p=f^p_i(\eta\rest n^p)$ for each $i\in s^p$ and for every $\eta\in2^{n^q}$;
	\item[(7)] $\gamma^q(\al)\rest n^p=\gamma^p(\al)$ for every $\al\in s^p$;
	\item[(8)] $\rho^q\rest{\dpower{v^p}2} =\rho^p$.
\end{enumerate}
It is easy to see (the details are given below) that the forcing $\poset$ introduces a countable family $\sett{f_n}{\ntr}$ of continuous functions on the Cantor set together with a function $\map \rho{\dpower{\omega_1}2}{\omega}$ and a one-to-one function $\map\gamma{\omega_1}{\cantor}$ such that $\gamma(\al)=f_{\rho(\al,\beta)}(\gamma(\beta))$ for every $\al<\beta<\omega_1$. 
We need to prove that $\poset$ does not collapse $\aleph_1$.

\begin{lm}\label{cece}
$\poset$ satisfies the countable chain condition.
\end{lm}

\begin{pf}
Fix a family $\Gee\subs\poset$ with $|\Gee|=\aleph_1$. Replacing $\Gee$ by an uncountable subfamily, we may assume that there exist $\ntr$, $s\in\fin\nat$ and $\Ef=\sett{f_i}{i\in s}\subs \nice(n)$ such that
$n^p=n$, $s^p=s$ and $\Ef^p=\Ef$ for every $p\in \Gee$. Further refining $\Gee$, we may assume that
\begin{enumerate}
	\item[(9)] $\setof{v^p}{p\in\Gee}$ forms a $\Delta$-system with root $a\subs\omega_1$.
	\item[(10)] For every $p,q\in\Gee$ the structures $\seq{v^p,\gamma^p,\rho^p,<}$ and $\seq{v^q,\gamma^q,\rho^q,<}$ are isomorphic, where $<$ is the linear order inherited from $\omega_1$. In other words, there exists an order preserving bijection $\map\phi{v^p}{v^q}$ such that $\gamma^p(\al)=\gamma^q(\phi(\al))$ and $\rho^p(\al,\beta)=\rho^q(\phi(\al),\phi(\beta))$ for every $\al,\beta\in v^p$.	
\end{enumerate}
For the remaining part of the proof we fix $p,q\in\Gee$ such that $\max(a)<\min(v^p\setminus a)$ and $\max(v^p)<\min(v^q\setminus a)$. Our aim is to construct $r\in\poset$ with $p\loe r$ and $q\loe r$.

Define $n^r=n+1$ and $v^r=v^p\cup v^q$. Note that by (10), $\rho^p$, $\rho^q$ coincide on $a=v^p\cap v^q$. 
Extend $\rho^p\cup\rho^q$ to a function $\map{\rho^r}{\dpower{v^r}{2}}\nat$ in such a way that $\rho^r$ restricted to the set
$$\sig = \dpower{v^r}{2} \setminus (\dpower{v^p}{2}\cup \dpower{v^q}{2})=\setof{\dn\al\beta}{\al\in v^p\setminus a, \;\; \beta\in v^q\setminus a}$$
is a bijection onto $t\subs\nat\setminus s$. Clearly, $\rho^r$ satisfies (2'), i.e. $\rho^(\al,\beta)\ne \rho^r(\al',\beta)$ whenever $\al<\al'<\beta$. Define $s^r=s\cup t$. Then $\map{\rho^r}{\dpower{v^r}{2}}{s^r}$.
Further, define
$$\gamma^r(\al)=
\begin{cases}
	\gamma^p(\al)\concat0 &\text{ if }\al\in v^p,\\
	\gamma^q(\al)\concat1 &\text{ if }\al\in v^q\setminus a.
\end{cases}$$
Observe that $\map{\gamma^r}{v^r}{2^{n^r}}$ is one-to-one.
It remains to define $\Ef^r=\sett{f^r_i}{i\in s^r}$.

If $i\in t$ then we define $f^r_i$ to be the constant function with value $\gamma^r(\al)$, where $\al\in v^p\setminus a$, $\beta\in v^q\setminus a$ are such that $i=\rho^r(\al,\beta)$. Note that $\al,\beta$ are uniquely determined, so there is no ambiguity here and $f^r_i$ satisfies (4). Finally, fix $i\in s$, $\eta\in2^n$, $\eps\in2$ and define
$$f^r_i(\eta\concat\eps)=
\begin{cases}
f_i(\eta)\concat\eps\quad& (\exists\;\al,\beta\in v^p\setminus a)\;\al<\beta \land i=\rho^p(\al,\beta) \land \eta=\gamma^p(\beta),\\
f_i(\eta)\concat0 &\text{ otherwise}.
\end{cases}
$$
By this way we have finished the definition of $r=\seq{n^r,s^r,v^r,\Ef^r,\gamma^r,\rho^r}$. 
In order to show that $r\in\poset$, we need to verify condition (4) only, since conditions (1)--(3) are rather clear.

For fix $\al<\beta$ in $v^r$ and let $\ell=\rho^r(\al,\beta)$. If $\ell\in t$ then $f^r_\ell$ is constantly equal to $\gamma^r(\al)$, therefore (4) holds in this case. So assume $\ell\in s$ and let $\eta=\gamma^r(\beta)\rest n$. 
We consider the following two cases.

\begin{case}
$\al\in v^q\setminus a$.
\end{case}

Notice that also $\beta\in v^q\setminus a$, because $\al<\beta$. By (10), there exist $\al',\beta'\in v^p$ such that $\rho^p(\al',\beta')=\rho^q(\al,\beta)=\ell$ and $\gamma^p(\beta')=\gamma^q(\beta)=\eta$. Thus the first possibility in the definition of $f^r_\ell$ occurs and we have $$f^r_\ell(\gamma^r(\beta))=f^r_\ell(\eta\concat1)=f_\ell(\eta)\concat1=\gamma^q(\al)\concat1 =\gamma^r(\al),$$ therefore (4) holds.

\begin{case}
$\al\in v^p$.
\end{case}

Now $\gamma^r(\al)=\gamma^p(\al)\concat0$ and either $\beta\in v^p$ or else $\al\in a$ and $\beta\in v^q$ (because $\ell\in s$ implies that either $\dn\al\beta\subs v^p$ or $\dn\al\beta\subs v^q$).
Observe that $f^r_\ell(\gamma^r(\beta))\rest n = \gamma^r(\al)\rest n$, by the definition of $f^r_\ell$ and by the fact that $p,q\in\poset$. Thus, the only possibility for the failure of (4) is that $f^r_\ell(\gamma^r(\beta))=\gamma^p(\al)\concat1$. Suppose this is the case. By the definition of $f^r_\ell$, we conclude that $\gamma^r(\beta)=\eta\concat1$ and in particular $\beta\in v^q\setminus a$ and $\al\in a$. 
Moreover, the first case in the definition of $f^r_\ell$ occurs, so there exist $\al'<\beta'$ in $v^p\setminus a$ such that $\ell=\rho^p(\al',\beta')$ and $\eta=\gamma^p(\beta')$.
Let $\map{\phi}{v^p}{v^q}$ be the bijection appearing in condition (10). In particular $\gamma^p(\beta)=\eta=\gamma^q(\phi(\beta'))$, therefore $\phi(\beta')=\beta$, because $\gamma^q$ is one-to-one.
Further, 
$$\rho^r(\phi(\al'),\beta)=\rho^q(\phi(\al'),\phi(\beta'))=\rho^p(\al',\beta')=\ell=\rho^r(\al,\beta).$$
Thus $\phi(\al')=\al$, because $\rho^r$ satisfies (2').
This leads to a contradiction, because $\al\in a$, $\al'\in v^p\setminus a$ and $\img\phi{v^p\setminus a}=v^q\setminus a$. Thus (4) holds.

We have proved that $r\in\poset$. Clearly $p\loe r$ and $q\loe r$.
\end{pf}

\begin{lm}\label{8a9}
Let $k\in\nat$ and $\xi\in\omega_1$. The sets
$$\Dee(k) = \setof{p\in\poset}{n^p\goe k\text{ and }k\in s^p},\qquad \Eee(\xi) = \setof{p\in\poset}{\xi\in v^p}.$$
are dense in $\poset$.
\end{lm}

\begin{pf}
Fix $p\in\poset$. Define $n^{q}=n^{p}+1$, $s^{q}=s^{p}\cup\sn k$, $v^{q}=v^{p}$, $\rho^q=\rho^p$, $\gamma^q(\eta)=\gamma^p(\eta)\concat0$ and $f^q_i(\eta\concat\eps)=f^p_i(\eta)\concat\eps$ for $i\in s^p$, $\eta\in 2^{n^p}$, $\eps\in2$. Finally, if $k\notin s^p$, let $f^q_k$ be any function from $\nice(n^p+1)$. By this way we have extended $p$ to a condition $q=\seq{n^q,s^q,v^q,\Ef^q,\gamma^q,\rho^q}\in \poset$ so that $n^q>n^p$, $k\in s^q$. Repeating this procedure finitely many times we obtain $r\goe p$ such that $n^r\goe k$ and $k\in s^r$. This shows that $\Dee(k)$ is dense in $\poset$.

In order to show the density of $\Eee(\xi)$ again fix $p\in\poset$ and assume $\xi\notin v^p$. Define $n^{q}=n^{p}+1$ and $v^{q}=v^{p}\cup\sn\xi$.
Let $\sig=\setof{\dn\xi\al}{\al\in v^p}$. Extend $\rho^p$ to a function $\map{\rho^q}{\dpower{v^q}2}\omega$ so that $\rho^q\rest\sig$ is one-to-one onto $t\subs\omega\setminus s^p$. Let $s^q=s^p\cup t$.
Further, define 
$\gamma^q(\al)=\gamma^p(\al)\concat0$ for $\al\in v^p$ and let $\gamma^q(\xi)$ be the constant one function in $2^{n^q}$. It remains to define $\Ef^q$.

Given $i\in s^p$, define $f^q_i(\eta\concat\eps)=f^p_i(\eta)\concat\eps$ for every $\eta\in 2^{n^p}$, $\eps\in2$. Fix $i\in t$ and let $\al\in v^p$ be such that $i=\rho^q(\al,\xi)$. If $\xi<\al$, define $f^q_i$ to be the constant function with value $\gamma^q(\xi)$. If $\al<\xi$, define $f^q_i$ to be the constant function with value $\gamma^q(\al)$. Observe that conditions (1) -- (4) are satisfied, therefore $q=\seq{n^q,s^q,v^q,\Ef^q,\gamma^q,\rho^q}\in \poset$. It is clear that $p\loe q$ and $q\in\Eee(\xi)$.
\end{pf}

\begin{lm}\label{kwadr4}
The poset $\poset$ forces a family $\Ef=\setof{f_n}{\ntr}$ of $1$-Lipschitz functions on the Cantor set $\cantor$ and an uncountable set $X\subs \cantor$ whose square is covered by $\Ef$.
\end{lm}

\begin{pf}
Let $G$ be a $\poset$-generic filter over a fixed ground model $\Ve$. Define functions $\map{f_k}{\cantor}{\cantor}$ ($k\in\nat$), $\map{\gamma}{\omega_1}{\cantor}$ and $\map{\rho}{\dpower{\omega_1}2}{\nat}$ by the following equations:
\begin{align*}
f_k(x)\rest n^p &= f_k^p(x\rest n^p),\\
\gamma(\al)\rest n^p &= \gamma^p(\al),\\
\rho(\al,\beta) &= \rho^p(\al,\beta),
\end{align*}
where $x\in \cantor$ and $p$ is any element of $G$ such that $\al,\beta\in v^p$ and $k\in s^p$. The fact that $G$ is a filter and the density of sets $\Dee(k)$ and $\Eee(\xi)$ (Lemma \ref{8a9}) imply that the above definitions are correct.
Let $X=\setof{\gamma(\xi)}{\xi<\omega_1}$. By the definition of $\poset$, the set $X\subs \cantor$ is uncountable, the functions $f_k$ are $1$-Lipschitz and for every $\al<\beta<\omega_1$ we have that $\gamma(\al)=f_{\rho(\al,\beta)}(\gamma(\beta))$.
It follows that $X^2\subs\sn{\id{\cantor}}\cup\bigcup_{\ntr}(f_n\cup f_n^{-1})$.
\end{pf}

\begin{pf}[Proof of Theorem \ref{wfojwqeoj}]
The above forcing shows that the existence of a family of $1$-Lipschitz functions on the Cantor set covering an uncountable square is consistent with ZFC. We are going to argue, using Keisler's completeness, that such a family actually exists in ZFC.

[...]
\end{pf}

\section{Final remarks}

It is natural to ask whether there exists an uncountable (necessarily scattered) compact space $K$ such that $K^2$ is covered by countably many graphs of continuous functions and their inverses. 
Below we show that the answer is negative.

\begin{tw}\label{ergjwej}
Let $K$ be a compact Hausdorff space and let $\ciag f$ be a family of continuous functions such that for each $\ntr$ the set $\dom(f_n)$ is closed in $K$ and $K\times K=\bigcup_{\ntr}(f_n\cup f_n^{-1})$. Then $|K|\loe\aleph_0$.
\end{tw}

\begin{pf}
By the Baire Category Theorem, a compact $K$ satisfying the above assertion must be scattered. Suppose the theorem is false and fix a counterexample $K$ of minimal Cantor-Bendixson rank $\lam$. Denote by $K^{(\al)}$ the $\al$-th derivative of $K$. Passing to a subspace, we may further assume that $K^{(\lam)}$ is a singleton, which we shall denote by $\infty$. 
Note that every closed set not containing $\infty$ is countable. Indeed, if $A\subs K$ is closed and $\infty\notin A$, then by compactness, $A\cap K^{(\gamma)}=\emptyset$ for some $\gamma<\lam$. Thus the Cantor-Bendixson rank of $A$ is $\loe\gamma$, therefore by the minimality of $\lam$, $A$ must be countable because it satisfies the above assertion.

Let $M=\setof{f_n(\infty)}{\ntr\text{ and }\infty\in\dom(f_n)}$ and choose $y\in K\setminus M$.
Let $$A = K\setminus (M\cup \setof{f_n(y)}{\ntr}).$$ Then $A$ is uncountable and for each $x\in A$ there exists $k\in\nat$ such that $y=f_k(x)$. Find $k\in\nat$ such that the set $B=\setof{x\in A}{y=f_k(x)}$ is uncountable. Note that $\infty\in\cl B$, because every closed set not containing $\infty$ is countable. Thus $\infty\in\dom(f_k)$ and, by continuity, $y=f_k(\infty)\in M$; a contradiction.
\end{pf}

By the above result, it is impossible to cover $\omega_1\times \omega_1$ by countably many functions which are continuous with respect to the order topology. Indeed, all these functions would be extendable onto the \v Cech-Stone compactification of $\omega_1$ which equals $\omega_1+1$ and therefore, adding one more function, we would obtain a countable family of continuous functions covering the square of $\omega_1+1$.

It is easy to see, using Sierpi\'nski's theorem, that the one point compactification of the discrete space of cardinality $\aleph_1$ can be covered by countably many {\em partial} continuous functions and their inverses. Thus, Theorem \ref{ergjwej} fails when we drop the assumption that $\dom(f_n)$ be closed.

\subsection*{Ackonwledgements}
The first author would like to thank S\l awomir Turek for useful comments and Piotr Zakrzewski for pointing out reference \cite{Zak}.

\end{document}